\documentclass[a4paper]{article}
\usepackage{graphicx} 
\usepackage{amsmath, amsthm, amsfonts, amssymb, amscd, bm}
\usepackage{verbatim}
\usepackage{enumerate}
\usepackage{url}

\newtheorem{theorem}{Theorem}[section]
\newtheorem{corollary}{Corollary}[theorem]
\newtheorem{lemma}{Lemma}[section]
\newtheorem{definition}{Definition}[section]

\theoremstyle{remark}
\newtheorem{remark}{Remark}[section]
\newtheorem{example}{Example}[section]

\DeclareMathOperator{\Tr}{Tr}

\begin{document}

\title{A General Theory of Finite State Backward Stochastic Difference Equations}
\author{Samuel N. Cohen \\ University of Adelaide \\ samuel.cohen@adelaide.edu.au 
\and Robert J. Elliott\thanks{Robert Elliott wishes to thank the Australian Research Council for
support.} \\ University of Calgary and University of Adelaide \\ relliott@ucalgary.ca}
\date{}

\maketitle
\begin{abstract}
By analogy with the theory of Backward Stochastic Differential Equations, we define Backward Stochastic Difference Equations on spaces related to discrete time, finite state processes. This paper considers these processes as constructions in their own right, not as approximations to the continuous case. We establish the existence and uniqueness of solutions under weaker assumptions than are needed in the continuous time setting, and also establish a comparison theorem for these solutions. The conditions of this theorem are shown to approximate those required in the continuous time setting. We also explore the relationship between the driver $F$ and the set of solutions; in particular, we determine under what conditions the driver is uniquely determined by the solution. Applications to the theory of nonlinear expectations are explored, including a representation result.
\end{abstract}

\section{Introduction}
The theory of Backward Stochastic Differential Equations (BSDEs) is an active area of research in both Mathematical Finance and Stochastic Control. Typically, one begins by defining processes $(Y, Z)$ through an equation of the form
\[Y_t - \int_{]t,T]} F(\omega, u, Y_{u-}, Z_u) du+ \int_{]t,T]} Z_u dM_u = Q.\]
Here $Q$ is a square-integrable terminal condition, $F$ a progressively measurable `driver' function, and $M$ an $N$-dimensional Brownian Motion, all defined on an appropriate filtered probability space. The `solutions' $(Y,Z)$ are required to be adapted to the forward filtration, and $Z$ is required to be predictable.

Recent work has also allowed the presence of jumps and the use of other underlying processes for $M$. However, these typically require a generalisation of the equation to include a martingale orthogonal to $M$, as a martingale representation theorem may not hold. See \cite{El1997a} for some general results. In \cite{Cohen2008}, we considered the situation where $M$ is the compensated jump martingale generated by a continuous-time, finite state Markov Chain and showed that solutions existed for equations of this type.

In this paper, we shall consider an analogous situation in discrete time. Such processes have been considered previously in \cite{Ma2002} and other works, particularly as numerical approximations to continuous time processes. In contrast to these authors, we approach discrete time BSDEs as entities in their own right, and do not significantly address their use as approximations and the related numerical methods. Because of this, we obtain considerably more general conditions under which solutions exist, and also establish fundamental results, including, for example, a comparison theorem. This helps provide a better understanding of the structure underlying our results, by removing the complexity of continuous time and the restrictions inherent in Brownian motion.

We begin by defining the discrete analogue of a continuous time BSDE, and giving conditions for existence and uniqueness. We then prove a comparison theorem, and consider the relationship between the driver and the set of solutions. We apply these results to obtain a theory of nonlinear expectations, and show that every nonlinear expectation in this context is indeed the solution to a discrete BSDE.

\section{Dynamics}

We shall consider an underlying discrete-time, finite state process $X$. Without loss of generality, this can be assumed to take values in the standard basis vectors of $\mathbb{R}^N$, where $N$ is the number of states of the process. That is, for each $t\in\{0,1,...\}$, 
\[X_t \in \{e_1, e_2,..., e_N\},\]
where $e_i = (0,0,..,0,1,0,...,0)^*\in \mathbb{R}^N$, and $[.]^*$ denotes vector transposition.

Let $(\Omega, \mathcal{F},\{\mathcal{F}_t\}_{0\leq t\leq T}, \mathbb{P})$ be a filtered probability space, where $\mathcal{F}_t$ is the completion of the sigma algebra generated by the process $X$ up to time $t$. Unless otherwise noted, we index all quantities by the first time $t$ such that they are $\mathcal{F}_t$-measurable.

$X$ can then be represented as:
\[X_t = E[X_t|\mathcal{F}_{t-1}] + M_t \in \mathbb{R}^N.\]
By definition, $M$ is the martingale difference process $M_t = X_t - E[X_t|\mathcal{F}_{t-1}]$. The central process considered in this paper is the solution $(Y,Z)$ of  a BSDE based on $M$, that is an equation of the form
\begin{equation}\label{eq:BSDEgen}
Y_t - \sum_{t\leq u< T} F(\omega, u, Y_u, Z_u) + \sum_{t\leq u< T} Z_u M_{u+1} = Q,
\end{equation}
where $T$ is a deterministic terminal time, $F$ an adapted map $F:\Omega\times \{0,...,T\} \times \mathbb{R}^K \times \mathbb{R}^{K\times N}\rightarrow \mathbb{R}^K$, and  $Q$ an $\mathbb{R}^K$ valued, $\mathcal{F}_T$-measurable terminal condition. 

We consider only solutions $(Y,Z)$ which are adapted to the filtration $\{\mathcal{F}_t\}$, that is, $(Y_t, Z_t)\in \mathbb{R}^K\times \mathbb{R}^{K\times N}$ is $\mathcal{F}_t$-measurable for all $t$. For simplicity, we shall also assume $Y_t,Z_t\in L^\infty(\mathcal{F}_t)$ for all $t$, $Q\in L^\infty(\mathcal{F}_T)$ and $F(\omega, t, Y_t, Z_t)\in L^\infty(\mathcal{F}_t)$ for all $t$ and all $Y_t, Z_t \in L^\infty(\mathcal{F}_t)$. Note, however, that as there are only finitely many possible paths for $X$ on $\{0,...,t\}$, and all our quantities are finite dimensional, it is clear that 
\[L^1(\mathcal{F}_t)=L^2(\mathcal{F}_t) = L^\infty(\mathcal{F}_t).\]

\begin{remark} As in \cite[Sec. 2]{Cohen2008b}, it is easy to extend this theory to the case where $T$ is an essentially bounded stopping time, by extending the domain of the driver $F$. For this reason, all the results obtained for deterministic $T$ considered here can be extended, after appropriate modification to the assumptions on $F$.
\end{remark}

\begin{definition}
For any integer $K$, we shall denote by $\|\cdot\|_M$ the seminorm on the space of adapted processes $Z$ in $\mathbb{R}^{K\times N}$, given by
\[\begin{split}
\|Z\|_M^2 &:= E\Tr\left[\sum_{0\leq u< T}Z_u\cdot E[M_{u+1} M_{u+1}^*|\mathcal{F}_{u}] \cdot Z_u^*\right]\\
&=\sum_{0\leq u< T}\Tr E[(Z_uM_{u+1})(Z_uM_{u+1})^*].\end{split}\]
\end{definition}

\begin{lemma} \label{lem:Msimequiv}
The following statements are equivalent:
\begin{enumerate}[(i)]
\item $\|Z^1- Z^2\|^2_M=0$.
\item $E\Tr[(Z^1_u-Z^2_u) M_{u+1} M_{u+1}^* (Z^1_u-Z^2_u)^*]=0$ for all $u\in\{0,...,T-1\}$.
\item $Z^1_uM_{u+1} = Z^2_u M_{u+1}$ $\mathbb{P}$-a.s. for all $u\in\{0,...,T-1\}$.
\item For all $t\in\{1,...,T\}$, $\mathbb{P}$-a.s.,
\[\sum_{0\leq u < t} Z^1_{u}M_{u+1} = \sum_{0\leq u < t} Z^2_{u}M_{u+1}\]
\end{enumerate}
In this case we shall write $Z^1\sim_M Z^2$.
\end{lemma}

\begin{proof} 
We know
\[\|Z^1-Z^2\|_M^2 =\sum_{0\leq u< T}\Tr E[((Z_u^1-Z_u^2)M_{u+1})((Z_u^1-Z_u^2)M_{u+1})^*].\]
For each $u$, the trace is then simply the sum of the expected squares of the components of $(Z_u^1-Z_u^2)M_{u+1}$, and so, each of the summed terms must be nonnegative. Hence the total sum is zero if and only if each term must be zero, that is,
 \[E\Tr[(Z^1_u-Z^2_u) M_{u+1} M_{u+1}^* (Z^1_u-Z^2_u)^*]=0 \text{ for all  }u\in\{0,...,T-1\}.\]
Therefore, (i) and (ii) are equivalent.

As a sum of squares, (ii) is true if and only if each term is zero, that is, 
\[E[(e_i^*(Z^1_u-Z^2_u) M_{u+1})^2]=0\]
for all basis vectors $e_i$. Considering all components at once, this is equivalent to $Z^1_uM_{u+1} = Z^2_u M_{u+1}$ $\mathbb{P}$-a.s. for all $u\in\{0,...,T-1\}$. Therefore, (ii) and (iii) are equivalent.

Taking a sum over the values $0\leq u<t$ in statement (iii) gives statement (iv). Taking the differences of statement (iv) at times $t$ and $t-1$ gives statement (iii).
\end{proof}

\begin{definition}

For two $\mathcal{F}_{t-1}$-measurable random variables $Z^1_{t-1}$ and $Z^2_{t-1}$, we shall write $Z^1_{t-1}\sim_{M_t} Z^2_{t-1}$ if  $Z^1_{t-1} M_t = Z^2_{t-1} M_t$, $\mathbb{P}$-a.s. 
\end{definition}

By Lemma \ref{lem:Msimequiv}, $Z^1\sim_M Z^2$ if and only if $Z^1_{t-1} \sim_{M_t} Z^2_{t-1}$ for all $t\in\{1,...,T\}$. It is straightforward to show that $\sim_M$ and $\sim_{M_t}$ are both equivalence relations.

\section{Existence and Uniqueness}
The first stage in any discussion of these processes is to establish under what conditions solutions of the BSDE (\ref{eq:BSDEgen}) exist and are unique.

Before proving the existence of these processes, we recall the following Martingale Representation Theorem (from \cite{Elliott1994}).
\begin{theorem}\label{thm:MartRep}
For any $\{\mathcal{F}_t\}$-adapted, $\mathbb{R}^K$ valued martingale $L$, there exists an adapted $\mathbb{R}^{K\times N}$ process $Z$ such that 
\[L_{t+1} = L_0 + \sum_{0\leq u < t} Z_u M_{u+1}.\]
This $Z$ process is unique up to equivalence $\sim_M$.
\end{theorem}
\begin{proof}
As $L$ is adapted, we can apply the Doob-Dynkin Lemma (see \cite[p174]{Shiryaev2000}) to show that, for each $t$, there exists some $\mathcal{F}_{t-1}$-measurable function $g_t:\mathbb{R}^N\rightarrow \mathbb{R}^K$ such that 
\[L_t = L_{t-1}+g_t(X_t).\]
Now $X_t$ can take $N$ possible values, associated with each of the basis vectors $e_i$ in $\mathbb{R}^N$. We can, therefore, create an $\mathcal{F}_{t-1}$-measurable $\mathbb{R}^{K\times N}$ matrix $Z_{t-1}$ with entries
\[Z_{t-1}=[g_t(e_1) |g_t(e_2)|...|g_t(e_N)]\]
which will satisfy
\[L_t = L_{t-1} + Z_{t-1} X_t.\]
$L_t$ is a martingale, hence $E[L_t|\mathcal{F}_{t-1}] = L_{t-1}$. Therefore, it must be the case that $Z_{t-1}E[X_t|\mathcal{F}_{t-1}] = 0$.
So \[Z_{t-1} M_t = Z_{t-1}[X_t-E[X_t|\mathcal{F}_{t-1}]] = Z_{t-1} X_t\]
giving
\[L_t = L_{t-1} + Z_{t-1}M_t.\]
Using this recursive formula as the basis for a telescoping sum gives the desired representation.

If we had two possible solutions, $Z^1$ and $Z^2$, then simple rearrangement gives
\[Z^1_{t-1}M_t=Z^2_{t-1}M_t\]
$\mathbb{P}$-a.s., for all $t$. Hence, by Lemma \ref{lem:Msimequiv}, $Z^1\sim_M Z^2$.
\end{proof}

\begin{corollary} \label{cor:meanzerorep}
For any $\mathbb{R}^K$ valued, $\mathcal{F}_t$-measurable random variable $W$ with $E[W|\mathcal{F}_{t-1}]=0$, there exists an $\mathcal{F}_{t-1}$-measurable random variable $Z_{t-1}$ such that, $\mathbb{P}$-a.s., \[W=Z_{t-1}M_{t}.\]
This variable is unique up to equivalence $\sim_{M_t}$.
\end{corollary}
\begin{proof}
We define a martingale $L$ by
\[L_s=I_{s<t}W.\]
It follows that $Z_s=0$ for $s<t-1$ and $s\geq t$, and the variable $Z_{t-1}$ is as desired.\end{proof}

\begin{theorem}\label{thm:BSDEExist}
Suppose $F$ is such that the following two assumptions hold:
\begin{enumerate}[(i)]
	\item For any $Y$, if $Z^1\sim_M Z^2$, then $F(\omega,t, Y_t, Z^1_t) = F(\omega, t, Y_t, Z^2_t)$ $\mathbb{P}$-a.s. for all $t$.
	\item For any $Z$, for all $t$, the map $Y_t\mapsto Y_t-F(\omega, t, Y_t, Z_u)$ is $\mathbb{P}$-a.s. a bijection from $\mathbb{R}^K\rightarrow \mathbb{R}^K$, up to equality $\mathbb{P}$-a.s.
\end{enumerate}
Then for any terminal condition $Q$ essentially bounded, $\mathcal{F}_T$-measurable, and with values in $\mathbb{R}^K$, the BSDE (\ref{eq:BSDEgen}) has an adapted solution $(Y, Z)$. Moreover, this solution is unique up to indistinguishability for $Y$ and equivalence $\sim_M$ for $Z$.
\end{theorem}
\begin{proof}
Clearly we can find an (adapted) solution $Y_T=Q$ at time $T$. We shall construct the solution for all $t$ using backward induction.

Suppose we have a solution at time $t+1$. In this case, equation (\ref{eq:BSDEgen}) can be simplified to
\begin{equation}\label{eq:OneStepBSDE}
Y_{t} - F(\omega, t, Y_t, Z_t) + Z_tM_{t+1} = Y_{t+1}.
\end{equation}
Taking a conditional expectation gives
\begin{equation}\label{eq:OneStepExpEq}
Y_t - F(\omega, t, Y_t, Z_t) = E[Y_{t+1}|\mathcal{F}_t],
\end{equation}
and so the martingale difference term must be $Y_{t+1}-E[Y_{t+1}|\mathcal{F}_t]$.
From Corollary \ref{cor:meanzerorep}, there is a $Z_t$ such that $Z_tM_{t+1}=Y_{t+1}-E[Y_{t+1}|\mathcal{F}_t]$. This $Z_t$ is unique up to equivalence $\sim_{M_{t+1}}$.

Using this $Z_t$, consider the equation
\[Y_t - F(\omega, t, Y_t, Z_t) = E[Y_{t+1}|\mathcal{F}_t].\]
This is now uniquely determined as an equation in $Y_t$. By assumption (ii), this equation has a unique solution $Y_t$, up to equality $\mathbb{P}$-a.s. The pair $(Y_t, Z_t)$ will solve the desired BSDE (i) at time $t$.

Backward induction starting with $T=t$ gives the desired result for all $t\in\{0,1,...T\}$. The solution $Y_t$ is unique up to equality $\mathbb{P}$-a.s. for all $t$, and hence $Y$ is unique up to indistinguishability (as we are working in discrete time). The solution $Z_t$ is unique up to equivalence $\sim_{M_{t+1}}$ for all $t$, and hence $Z$ is unique up to equivalence $\sim_M$.
\end{proof}

\begin{remark}
Note that, unlike in the continuous time case, we have not assumed any continuity conditions on $F$, Lipschitz or otherwise. In particular, the assumptions on $F$ as a function of $Z$ are very weak, essentially demanding only that $F$ does not distinguish ``equivalent'' strategies $Z$.
\end{remark}
\begin{corollary} \label{cor:BSDEExistAssnNec}
In general, both assumptions of Theorem \ref{thm:BSDEExist} are necessary for a unique solution to exist. That is, if either assumption fails for some $Y$ or $Z$, then there will exist terminal conditions $Q$ such that either no solutions or multiple solutions of the BSDE (\ref{eq:BSDEgen}) exist (on some set $t\in\{s,...T\}$).
\end{corollary}
\begin{proof}
{\bf Necessity of (i)}: Suppose that for some $Y$, some $Z$, $\tilde Z$ with $Z\sim_M \tilde Z$, \[F(\omega, t, Y_t, Z_t)\neq F(\omega, t, Y_t, \tilde Z_t)\] for some $t$ with some positive probability. Note $Y_t$ satisfies 
\[Y_t - F(\omega, t, Y_t, Z_t)= E[Y_{t+1}].\]

As $Z \sim_M \tilde Z$, we know that $Z_t M_{t+1} = Z_t M_{t+1}$ $\mathbb{P}$-a.s. Let $\tilde Y_t$ be a solution to the equation
\[\tilde Y_t - F(\omega, t, \tilde Y_t, \tilde Z_t)= E[Y_{t+1}].\]
We have assumed $F(\omega, t, Y_t, Z_t)\neq F(\omega, t, Y_t, \tilde Z_t)$ with positive probability, and it follows that $\tilde Y_t\neq Y_t$ with positive probability. Nevertheless, we have 
\[Y_t - F(\omega, t, Y_t, Z_t) + Z_t M_{t+1}=Y_{t+1}=\tilde Y_t - F(\omega, t, \tilde Y_t, \tilde Z_t) + \tilde Z_t M_{t+1}.\]
As there are two distinguishable solutions $Y_t$ and $\tilde Y_t$ to (\ref{eq:OneStepBSDE}), the one-step version of the BSDE (\ref{eq:BSDEgen}), the solution of (\ref{eq:BSDEgen}) is not unique.

{\bf Necessity of (ii) -- Injectivity:} Suppose, for some $t$ and some $Z$, the map $Y_t\mapsto Y_t-F(\omega, t, Y_t, Z_u)$ is not injective with positive probability. Then there exist two distinguishable values $Y_t$, $\tilde Y_t$ such that $Y_t -F(\omega, t, Y_t, Z_t) = \tilde Y_t -F(\omega, t, \tilde Y_t, Z_t)$ $\mathbb{P}$-a.s., and these will both solve the desired equation (\ref{eq:OneStepBSDE})
\[Y_t - F(\omega, t, Y_t, Z_t) + Z_t M_{t+1}=Y_{t+1}=\tilde Y_t - F(\omega, t, \tilde Y_t, \tilde Z_t) + \tilde Z_t M_{t+1},\]
and hence will lead to distinguishable solutions $Y$ and $\tilde Y$ of the BSDE (\ref{eq:BSDEgen}).

{\bf Necessity of (ii) -- Surjectivity:} Suppose that for some $t$, some $Z$ the map $Y_t\mapsto Y_t-F(\omega, t, Y_t, Z_u)$ is not surjective with positive probability. Then there exists a value $q\in\mathbb{R}^K$ such that the equation $q = Y_t -F(\omega, t, Y_t, Z_t)$ has no solution $Y_t$ with positive probability. Hence, for any terminal condition $Q$ with $E[Y_{t+1}]=q$, there exists no $\mathcal{F}_t$-measurable $Y_t$ which will $\mathbb{P}$-a.s. satisfy (\ref{eq:OneStepExpEq}), and hence no adapted solution $Y$ to the BSDE (\ref{eq:BSDEgen}). (Note that such a terminal condition is clear if $T=t+1$ and $Q=Y_{t+1}$.)
\end{proof}

\begin{remark}
When these assumptions fail, it is possible for there to be a pathological situation where, for some $r<s$, multiple solutions ($Y_s$) of the BSDE (\ref{eq:BSDEgen}) exist for $t=s$, but for only one of these values of $Y_s$ do solutions exist to the equation
\[Y_r - \sum_{r\leq u< s} F(\omega, u, Y_u, Z_u) + \sum_{r\leq u< s} Z_u M_{u+1} = Y_s.\]
In this case we would have to say that there is a unique solution for $Y$ on $\{r,...,T\}$, but there are multiple solutions for $Y$ on $\{s,...T\}$.

The conditions of Theorem \ref{thm:BSDEExist} are therefore necessary and sufficient for unique solutions to exist for {\em all} values of $t$.
\end{remark}

\begin{remark} \label{rem:lipschitz}
In \cite{Ma2002}, a similar existence theorem for a solution $(Y,Z)$ is proven using the assumption that $F$ is Lipschitz, and a fixed point argument is used. This was motivated by the fact that the processes in \cite{Ma2002} are approximations of continuous time processes, and, therefore, the driver is actually of the form $\tilde F/n$, where $\tilde F$ is the (Lipschitz) continuous-time driver, and $n$ is the number of time steps between $0$ and $1$ used in the discrete approximation. Considering the scalar case, this implies that, if, $\mathbb{P}$-a.s.,
\[|\tilde F(\omega, t, Y^1, Z) - \tilde F(\omega, t, Y^2, Z)| \leq c |Y^1-Y^2|\]
then, for any $n>c$, the map 
\[Y \mapsto Y-\tilde F(\omega, t, Y, Z)/n\]
is $\mathbb{P}$-a.s. a strictly increasing function, and hence, is bijective. If $\tilde F$ is also appropriately Lipschitz continuous in $Z$, (see, for example, the conditions of \cite[Eq. 5.8]{El1997a} or \cite[Section 6]{Cohen2008}), assumption (i) of Theorem \ref{thm:BSDEExist} will also be satisified.  Hence, by Theorem \ref{thm:BSDEExist}, a unique solution will exist.

The assumptions of Theorem \ref{thm:BSDEExist}, on the other hand, are considerably more general, as they allow discontinuities in $F$, provided the stated bijective property holds. In particular, no continuity in $F$ as a function of $Z$ is assumed. Corollary \ref{cor:BSDEExistAssnNec} states the assumptions of Theorem  \ref{thm:BSDEExist} are also the most general conditions under which a unique solution of (\ref{eq:BSDEgen}) will exist for all $Q$.
\end{remark}

The following lemma will prove useful later.

\begin{lemma} \label{lem:bijective}
Let $(Y,Z)$ be the solution to a discrete BSDE satisfying Theorem \ref{thm:BSDEExist}. Then, for each $t$, there exists a bijection, (up to equality $\mathbb{P}$-a.s. and equivalence $\sim_{M_{t+1}}$), between the pair $(Y_t, Z_t)$ and the value $Y_{t+1}$.

Furthermore, this implies that there exists a bijection, (up to equality $\mathbb{P}$-a.s. and equivalence $\sim_{M}$), between $(Y_0, \{Z_s\}_{s<t})$ and $Y_t$, and that there exists a bijection, (up to equality $\mathbb{P}$-a.s. and equivalence $\sim_{M}$), between $(Y_t, \{Z_s\}_{t<T})$ and $Y_T$.
\end{lemma}
\begin{proof}
For each pair $(Y_t, Z_t)$, there exists a unique $Y_{t+1}$ given by 
\begin{equation} \label{eq:OneStepBij}
Y_{t+1} = Y_t- F(\omega, t, Y_t, Z_t)+Z_t M_{t+1}.
\end{equation}

For any $Y_{t+1}\in L^\infty(\mathcal{F}_{t+1})$, we wish to show there is a unique pair $(Y_t, Z_t)$ such that (\ref{eq:OneStepBij}) is true. It follows from (\ref{eq:OneStepBij}), and the fact $F(\omega, t, Y_t, Z_t)$ is $\mathcal{F}_t$-measurable, that $Z_tM_{t+1}=Y_{t+1}-E[Y_{t+1}|\mathcal{F}_t]$. By Corollary \ref{cor:meanzerorep}, there is a unique such $Z_t$, up to equivalence $\sim_{M_{t+1}}$. Finally, it follows from (\ref{eq:OneStepBij}) that $Y_t- F(\omega, t, Y_t, Z_t)=E[Y_{t+1}|\mathcal{F}_t]$. We have assumed the left hand side of this equation is bijective as a function of $Y_t$ and therefore there is a unique $(Y_t, Z_t)$ which satisfies (\ref{eq:OneStepBij}), for a given $Y_{t+1}$. Therefore, there is a bijective relationship between the pair $(Y_t, Z_t)$ and $Y_{t+1}$.

As there exists a bijection between the pair $(Y_t, Z_t)$ and $Y_{t+1}$ for all $t$, the remaining statements clearly follow by induction.
\end{proof}

\section{A Comparison Theorem}
A key result in the study of BSDEs is the `Comparison Theorem', first obtained in \cite{Peng1992}. As is shown in \cite{Cohen2008b}, for Markov chain driven BSDEs, there are conditions under which a comparison theorem remains true when the underlying process is not a Brownian motion. We now present a comparison theorem for discrete time BSDEs. For ease of notation we make the following definition:
\begin{definition}
Let $\mathbb{J}_t$ denote the $\mathcal{F}_t$-measurable set of indices of possible values of $X_{t+1}$, given $\mathcal{F}_{t}$. That is,
\[\mathbb{J}_t :=\{i:\mathbb{P}(X_{t+1}=e_i | \mathcal{F}_t)>0\}.\]
\end{definition}

In the following, an inequality on a vector quantity is to hold componentwise. 

\begin{theorem}[Comparison Theorem]\label{thm:CompThm}
Consider two discrete time BSDEs as in (\ref{eq:BSDEgen}), corresponding to coefficients $F^1, F^2$ and terminal values $Q^1, Q^2$. Suppose the conditions of Theorem \ref{thm:BSDEExist} are satisfied for both equations, let $(Y^1, Z^1)$ and $(Y^2, Z^2)$ be the associated solutions. Suppose the following conditions hold:
\begin{enumerate}[(i)]
	\item $Q^1\geq Q^2$ $\mathbb{P}$-a.s.
	\item $\mathbb{P}$-a.s.,  for all times $t$, \[F^1(\omega, t, Y_t^2, Z_t^2) \geq F^2(\omega, t, Y_t^2, Z_t^2).\]
	\item $\mathbb{P}$-a.s., for all $t$, for all $i$, the $i$th component of $F^1$, given by $e_i^*F^1$, satisfies
	\[e_i^*F^1(\omega, t, Y_t^2, Z_t^1) - e_i^*F^1(\omega, t, Y_t^2, Z_t^2)\geq\min_{j\in \mathbb{J}_t}\{e_i^*(Z^1_t-Z^2_t)(e_j-E[X_{t+1}|\mathcal{F}_{t}])\}.\]
	\item $\mathbb{P}$-a.s., for all $t$, if
	\[Y^1_t -F^1(\omega, t, Y_t^1, Z_t^1) \geq Y^2_t-F^1(\omega, t, Y_t^2, Z_t^1)\]
	then $Y^1_t\geq Y^2_t$.
\end{enumerate}
 It is then true that $Y^1 \geq Y^2$ $\mathbb{P}$-a.s. 
\end{theorem}
\begin{proof}
We shall establish this theorem using backward induction. For $t=T$ it is clear that $Y_t^1-Y_t^2=Q^1-Q^2\geq 0$ $\mathbb{P}$-a.s. as desired. Recall 
\[Y^i_t - \sum_{t\leq u<T} F^i(\omega, u, Y^i_u, Z^i_u)+\sum_{t\leq u<T} Z_u^iM_{u+1} = Q^i\]
for $i=1,2$. Then taking the one step equation, as in (\ref{eq:OneStepBSDE}) we have
\[Y^i_{t} - F^i(\omega, t, Y^i_t, Z^i_t) + Z^i_tM^i_{t+1} = Y^i_{t+1}\]
for all $0\leq t<T.$

For a given $t$, suppose we know $Y_{t+1}^1-Y_{t+1}^2\geq 0$ $\mathbb{P}$-a.s. Then, omitting the $\omega$ and $t$ arguments of $F^1$ and $F^2$, 
\[Y^1_t-Y^2_t - F^1(Y_t^1, Z_t^1) + F^2(Y_t^2, Z_t^2) + (Z_t^1-Z_t^2)M_{t+1} = Y^1_{t+1}-Y^2_{t+1}\geq 0.\]
Recalling $M_{t+1}=X_{t+1}-E[X_{t+1}|\mathcal{F}_t]$ and that $X_{t+1}$ takes values from the basis vectors $e_i$, we see that, for each component $e_i$,  $\mathbb{P}$-a.s.,
\[e_i^*(Y^1_t-Y^2_t) \geq e_i(F^1(Y_t^1, Z_t^1) - F^2(Y_t^2, Z_t^2)) - \min_{j\in \mathbb{J}_t}\{e_i^*(Z_t^1-Z_t^2)(e_j-E[X_{t+1}|\mathcal{F}_t])\}.\]
Hence, again $\mathbb{P}$-a.s., Assumptions (ii) and (iii) imply
\begin{equation}\label{eq:splitup}
\begin{split}
&e_i^*(Y^1_t-Y^2_t -F^1(Y_t^1, Z_t^1) + F^1(Y_t^2, Z_t^1))\\
&\geq e_i^*(F^1(Y_t^2, Z_t^2) - F^2(Y_t^2, Z_t^2))\\&\qquad+e_i^*F^1(Y_t^2, Z_t^1) - e_i^*F^1(Y_t^2, Z_t^2) - \min_{j\in \mathbb{J}_t}\{e_i^*(Z_t^1-Z_t^2)(e_j-E[X_{t+1}|\mathcal{F}_t])\}\\
&\geq 0.\end{split}
\end{equation}
That is, the inequality being taken componentwise, \[Y^1_t-Y^2_t -F^1(Y_t^1, Z_t^1) + F^1(Y_t^2, Z_t^1)\geq 0,\] and hence, by Assumption (iv),
\[Y^1_t\geq Y^2_t\]
$\mathbb{P}$-a.s. as desired. The general statement follows by backward induction.
\end{proof}

\begin{corollary} \label{cor:StrictComp}
Suppose Theorem \ref{thm:CompThm} holds and the inequality in Assumption (iii) is strict, that is, for all $i$, 
	\[e_i^*F^1(\omega, t, Y_t^2, Z_t^1) - e_i^*F^1(\omega, t, Y_t^2, Z_t^2)>\min_{j\in\mathbb{J}_t}\{e_i^*(Z^1_t-Z^2_t)(e_j-E[X_{t+1}|\mathcal{F}_{t}])\}\]
	unless $e_i^*Z_t^1 \sim_{M_{t+1}} e_i^*Z_t^2$.
	
 Then this comparison is strict, that is, if on some $A\in \mathcal{F}_t$ we have $Y^1_t = Y^2_t$ $\mathbb{P}$-a.s. on $A$, then $Q^1=Q^2$ $\mathbb{P}$-a.s. on $A$, and for all $s\in\{t,...,T\}$, $\mathbb{P}$-a.s. on $A$, $F^1(\omega, s, Y_s^2, Z_s^2) = F^2(\omega, s, Y_s^2, Z_s^2)$, $Z_s^1\sim_{M_{s+1}}Z_s^2$ and $Y^1_s = Y^2_s$.
\end{corollary}
\begin{proof}
Throughout this proof we will omit the $\omega$ and $t$ arguments of $F^1$ and $F^2$, and all (in-)equalities are assumed to hold $\mathbb{P}$-a.s. on $A$.

 In this case, for a given $t$, by the same argument as used to show (\ref{eq:splitup}), we can establish the strict inequality, for each $i$,
\begin{equation}\label{eq:splitup2}
\begin{split}
&e_i^*(Y^1_t-Y^2_t -F^1(Y_t^1, Z_t^1) + F^1(Y_t^2, Z_t^1))\\
&= e_i^*(Y^1_{t+1}-Y^2_{t+1})+e_i^*(F^1(Y_t^2, Z_t^2) - F^2(Y_t^2, Z_t^2))\\&\qquad+e_i^*F^1(Y_t^2, Z_t^1) - e_i^*F^1(Y_t^2, Z_t^2) - e_i^*(Z_t^1-Z_t^2)M_{t+1}\\
&\geq e_i^*(F^1(Y_t^2, Z_t^2) - F^2(Y_t^2, Z_t^2))\\&\qquad+e_i^*F^1(Y_t^2, Z_t^1) - e_i^*F^1(Y_t^2, Z_t^2) - \min_{j\in \mathbb{J}_t}\{e_i^*(Z_t^1-Z_t^2)(e_j-E[X_{t+1}|\mathcal{F}_t])\}\\
&> 0,\end{split}
\end{equation}
unless $Z_t^1\sim_{M_{t+1}} Z_t^2$. Hence, if $Y^1_t=Y^2_t$, the first term of this inequality is zero, which is only the case if $Z_t^1\sim_{M_{t+1}} Z_t^2$. 

If $Z_t^1\sim_{M_{t+1}} Z_t^2$, we know that 
\[F^1(Y_t^2, Z_t^1) - F^1(Y_t^2, Z_t^2) - (Z_t^1-Z_t^2)M_{t+1}=0\]
and so, taking a conditional expectation through (\ref{eq:splitup2}), and combining the results for each $i$,
\[\begin{split}
0&\geq -E[Y^1_{t+1}-Y^2_{t+1}|\mathcal{F}_t]\\
&=Y^1_t-Y^2_t -F^1(Y_t^1, Z_t^1) + F^1(Y_t^2, Z_t^1)\\
&\geq F^1(Y_t^2, Z_t^2) - F^2(Y_t^2, Z_t^2)\end{split}\]
the inequality being taken componentwise. As the final term is nonnegative, this can only be satisfied if 
\[F^1(Y_t^2, Z_t^2) = F^2(Y^2_t, Z_t^2).\] 
and $Y^1_{t+1} = Y^2_{t+1}$. The result follows by forward induction.
\end{proof}

\begin{corollary} \label{cor:CompThmAssnSwitch}
Theorem \ref{thm:CompThm} remains true if we replace Assumptions (iii) and (iv) by 

\begin{enumerate}[(i')]
	\setcounter{enumi}{2}
	\item $\mathbb{P}$-a.s., for all $t$, for all $i$, $e_i^*F^1$ satisfies
	\[e_i^*F^1(\omega, t, Y_t^1, Z_t^1) - e_i^*F^1(\omega, t, Y_t^1, Z_t^2)\geq\min_{j\in \mathbb{J}_t}\{e_i^*(Z^1_t-Z^2_t)(e_j-E[X_{t+1}|\mathcal{F}_{t}])\}.\]
	\item $\mathbb{P}$-a.s., for all $t$, if
	\[Y^1_t -F^1(\omega, t, Y_t^1, Z_t^2) \geq Y^2_t-F^1(\omega, t, Y_t^2, Z_t^2)\]
	then $Y^1_t\geq Y^2_t$.
\end{enumerate}
or with a mixture of (iii-iv) and (iii'-iv') for different times $t$, provided one of these pairs holds for each $t$.
Similarly for the strict comparison of Corollary \ref{cor:StrictComp}.
\end{corollary}
\begin{proof}
The proof is essentially the same, except that, when using (iii') and (iv'), the decomposition in (\ref{eq:splitup}) becomes
\[
\begin{split}
&e_i^*(Y^1_t-Y^2_t) -e_i^*(F^1(Y_t^1, Z_t^2) + F^1(Y_t^2, Z_t^2))\\
&\geq e_i^*(F^1(Y_t^2, Z_t^2) - F^2(Y_t^2, Z_t^2))\\&\qquad+e_i^*F^1(Y_t^1, Z_t^1) - e_i^*F^1(Y_t^1, Z_t^2) - \min_{j\in \mathbb{J}_t}\{e_i^*(Z_t^1-Z_t^2)(e_j-E[X_{t+1}|\mathcal{F}_t])\}\\
&\geq 0,\end{split}
\]
and similarly in (\ref{eq:splitup2}). The remainder of the proof follows as before.
\end{proof}

\begin{remark}
Assumption (iii) of Theorem \ref{thm:CompThm} is closely related to the assumption used in \cite[Thm 4.2]{Cohen2008b} for continuous time BSDEs driven by Markov Chains. It essentially ensures that there will always be a value of $M_{t+1}$, that occurs with positive probability, such that \[-F^1(\omega, t, Y_t^2, Z_t^1) + F^1(\omega, t, Y_t^2, Z_t^2)+(Z^1_t-Z^2_t)M_{t+1}\]
is negative.
\end{remark}

\begin{remark}
Unlike the corresponding theorem in \cite{Peng1992}, Theorem \ref{thm:CompThm} holds for both scalar, ($K=1$), and vector, ($K>1$), valued $Y$. The key assumption in this regard is Assumption (iv), which is a signficantly more restrictive assumption in the vector case.

In the scalar case, Assumption (iv) is simply that the map $\phi:Y\mapsto Y-F^1(\omega, t, Y, Z^1_t)$ is $\mathbb{P}$-a.s. strictly increasing. If $F^1$ is differentiable (and hence continuous) with respect to $Y$ and has derivative below $1$ at any point, then, as we have assumed $\phi$ is a bijection, we can see that $\phi$ must be strictly increasing.

In the vector case, Assumption (iv) is nontrivial, even when approximating a continuous time process. Examples of relevant continuous time assumptions can be seen in the context of Markov chain drivers in \cite[Thm 5.3, 5.7]{Cohen2008b} and the associated discussion. A vector counterexample to the comparison theorem when Assumption (iv) fails is given below.
\end{remark}

\begin{example}
For approximating a univariate Brownian motion for a scalar BSDE we can see the following:

A simple approximation for a univariate Brownian motion is for each $X_{t}$ to be in one of two states, with equal probability, independently of the past. We then consider $Z$ such that the two values of $Z_tM_{t+1}$ are equal in magnitude and opposite in sign, and are of the order of $1/\sqrt{n}$, that is,
\[Z_t = z_t[1/\sqrt{n}, -1/\sqrt{n}]\]
for some real valued, continuous time, predictable process $z$. This binomial random walk model is considered, though with different notation, in \cite{Ma2002}. With this choice of $Z$, it is possible to show that 
\[\sum_{t\leq u<T}Z_uM_{u+1} \rightarrow \int_{]t,T]} z_{t}dW_t\]
in some sense, as $n\rightarrow \infty$, for $W$ a standard Brownian motion.

 When we are using a discrete process to approximate a continuous one, for $i=1,2$, $F^i$ will often be of the form \[F^i(\omega,t, Y_t, Z_t) = \tilde F^i(\omega, t, Y_t, z_t)/n,\] for $\tilde F^i$ the Lipschitz continuous driver in the continuous time equation, (as in \cite{Ma2002}). As noted in Remark \ref{rem:lipschitz} above, if $\tilde F$ has Lipschitz constant $c$ then Assumption (iv)of Theorem \ref{thm:CompThm} will be satisfied as soon as $n>c$, that is, provided we are approximating on a fine enough grid in the time dimension.

Then, $F^1=\tilde F^1 /n$ as above and for $n$ sufficiently large we will have 
\[\begin{split}
F^1(\omega, t, Y^1_t, Z^1_t)-F^1(\omega, t, Y^1_t, Z^2_t)&=(\tilde F^1(\omega, t, Y^1_t, z^1_t)-\tilde F^2(\omega, t, Y^1_t, z^2_t))/n \\
&>-|z^1_t-z^2_t|/\sqrt{n}\\&=\min_i\{(Z^1_t-Z^2_t)(e_i-E[X_{t+1}|\mathcal{F}_t]\},\end{split}\]
 and so it is clear that Assumption (iii) of Theorem \ref{thm:CompThm} will also be satisfied.
\end{example}

We now present counterexamples to Theorem \ref{thm:CompThm} when one of Assumptions (iii) or (iv) fails.
\begin{example}
Consider a pair of scalar BSDEs satisfying Theorem \ref{thm:BSDEExist}, with terminal conditions $Q^1=Q^2$, and driver $F^1=F^2=F$. Assume that $Y_0\mapsto Y_0-F(\omega, 0, Y_0, Z_0^1)$ is a strictly increasing function of $Y$. For simplicity, suppose the terminal time is $T=1$. Suppose Assumption (iii) of Theorem \ref{thm:CompThm} does not apply, in particular, that 
\[\begin{split}0<-F(\omega, 0, Y_0^2, Z_0^1) + F(\omega, 0, Y_0^2, Z_0^2)+\min_{i\in\mathbb{J}_0}\{[Z^1_0-Z^2_0]^*(e_i-E[X_1|\mathcal{F}_{0}])\}. \end{split}\]

Then we have, $\mathbb{P}$-a.s.,
\[\begin{split}0&=Y_1^1-Y^2_1\\
&=Y^1_0-Y^2_0 -F(\omega, 0, Y_0^1, Z_0^1) + F(\omega, 0, Y_0^2, Z_0^2) + (Z_0^1-Z_0^2)M_1\\
&>Y^1_0-Y^2_0 - F(\omega, 0, Y_0^1, Z_0^1) + F(\omega, 0, Y_0^2, Z_0^1)
\end{split}\]
and hence 
\[Y^2_0 -  F(\omega, 0, Y_0^2, Z_0^1) > Y^1_0- F(\omega, 0, Y_0^1, Z_0^1).\]
As the map $Y_0\mapsto Y_0- F(\omega, 0,  Y_0, Z_0^1)$ is assumed to be strictly increasing, this shows that $Y^2_0>Y^1_0$, contradicting the conclusion of Theorem \ref{thm:CompThm}.
\end{example}

\begin{example}
Consider a pair of vector valued BSDEs with $K=2$. Again assume $T=1$. For any $\mathbb{R}^2$ valued function $f$, let $F^1=F^2=F$ with
\[F(\omega, 0, Y_0, Z_0) = f(\omega) + \left[\begin{array}{cc}0&-1\\0&0\end{array}\right]Y_0.\]
Note that Assumptions (ii) and (iii) of Theorem \ref{thm:CompThm} are trivially satisfied.

Suppose $Q^1-Q^2$ satisfies 
\[E[Q^1-Q^2] = \left[\begin{array}{c} 0\\5\end{array}\right].\]
Then, taking an expectation of the difference of the BSDEs it is easy to show that $Y^1_0-Y^2_0$ satisfies
\[Y^1_0-Y^2_0 - \left[\begin{array}{cc}0&-1\\0&0\end{array}\right](Y^1_0-Y^2_0) = \left[\begin{array}{c} 0\\5\end{array}\right],\]
which implies
\[Y^1_0-Y^2_0 = \left[\begin{array}{c} -5\\5\end{array}\right].\]

Therefore, it is clear that $Y^1_0$ is not greater than $Y^2_0$ componentwise, contradicting the conclusion of Theorem \ref{thm:CompThm}.

Changing the matrix in the definition of $F$ in this example can also lead to other behaviour, for example with
\[F(\omega, 0, Y_0, Z_0) = f(\omega) + \left[\begin{array}{cc}-2&0\\0&-2\end{array}\right]Y_0\]
and $E[Q^1-Q^2] = [1,1]^*$ we find 
\[Y^1_0-Y^2_0 = \left[\begin{array}{c} -1\\-1\end{array}\right],\]
a complete reversal of the sign of $Y^1-Y^2$ in every component.
\end{example}

\section{Observing the driver}
Usually, the function $F$ is given, and our task is to obtain solutions to the (discrete or continuous time) BSDE. It may be of interest in applications to consider the reverse problem, that is to see whether, given the solutions to the BSDE, we can determine the values of the function $F$. A short comment on doing this is in \cite{Peng2005}, with other results in \cite{Briand2000}. However, this is only in continous time with an underlying Brownian motion, and, as might be expected, requires various limiting arguments. We here show that, in the discrete time context considered here, we can explicitly determine the function $F$.

\begin{definition}
Let $F$ be a driver for a discrete time BSDE (\ref{eq:BSDEgen}) satisfying the conditions of Theorem \ref{thm:BSDEExist}. We define $\mathfrak{O}_t$, the ``one step values'' under $F$,  to be the set of ordered pairs
\[\mathfrak{O}_t:=\{(Y_t, Y_{t+1})\},\] where $Y_t$ is the time $t$ solution to the BSDE with driver $F$ and terminal condition $Y_{t+1}$ at time $t+1\leq T$.
\end{definition}

Note that these values are not arbitrary, but are known to have arisen from some BSDE. It immediately apparent that this implies each terminal value $Y_{t+1}\in L^\infty(\mathcal{F}_{t+1})$ appears once and only once in $\mathfrak{O}_t$. In Theorem \ref{thm:phibijectBSDE}, we present a necessary and sufficient condition on $\mathfrak{O}_t$ for this to be true.

\begin{theorem} \label{thm:OdetF}
For a given $\mathfrak{O}_t$, there is a unique function $F(\omega, t, \cdot, \cdot)$ associated with $\mathfrak{O}_t$, that is, where $F(\omega, t, Y_t, Z_t)$ is the value at time $t$, for a given $Y_t$ and $Z_t$, of the driver $F$ generating $\mathfrak{O}_t$.
\end{theorem}

\begin{proof}
We know that $\mathfrak{O}_t$ comes from some BSDE and, therefore, there exists a function $F$ which generates it.

Consider $Y_{t+1}\in L^\infty(\mathcal{F}_{t+1})$. For each such $Y_{t+1}$, there is a unique associated value
\[(Y_t, Y_{t+1})\in\mathfrak{O}_t.\]
We know $Z_tM_{t+1}=Y_{t+1}-E[Y_{t+1}|\mathcal{F}_t]$ uniquely defines a value $Z_t$ (up to equivalence $\sim_{M_{t+1}}$), and that this value of $Z_t$ will be the time $t$ value of the solution to the BSDE (\ref{eq:BSDEgen}). By Lemma \ref{lem:bijective} there exists a $Y_{t+1}$ associated with every pair $(Y_t, Z_t)$ and therefore, 
\[F(\omega, t, Y_t, Z_t) = Y_{t+1}-E[Y_{t+1}|\mathcal{F}_t]\]
uniquely determines $F$ at time $t$.
\end{proof}

\begin{remark}
In a financial context, this implies that, assuming prices are generated by a BSDE, if a price $Y_t$ is given for every asset $Y_{t+1}$ when it is sold at time $t+1$, then the full behaviour of the BSDE can be determined.
\end{remark}

The following theorem acts as a converse to Lemma \ref{lem:bijective}, and gives necessary and sufficient conditions for $\mathfrak{O}_t$ to arise from a BSDE.

\begin{theorem} \label{thm:phibijectBSDE}
For every $Y_{t+1}\in L^\infty(\mathcal{F}_{t+1})$ fix a pair $(Y_t, Y_{t+1})$. Then the following statments are equivalent:
\begin{enumerate}
\item There is a map $\phi$ (unique up to equality $\mathbb{P}$-a.s.)
\[\phi: L^\infty(\mathcal{F}_{t+1})\rightarrow L^\infty(\mathcal{F}_{t}), Y_{t+1}\mapsto Y_t\] 
such that, for all $Z_t$,
\[k\mapsto \phi_{Z_t}(k):=\phi(k+Z_tM_{t+1})\]
is $\mathbb{P}$-a.s. a bijection in $k\in L^\infty(\mathcal{F}_t)$, and is invariant under equivalence $\sim_{M_{t+1}}$ for $Z_t$.
\item There exists a driver $F$ (unique up to equality $\mathbb{P}$-a.s.), satisfying the conditions of Theorem \ref{thm:BSDEExist}, such that 
\[Y_{t+1} = Y_t - F(\omega, t, Y_t, Z_t) +Z_t M_{t+1}.\]
\end{enumerate}
\end{theorem}
\begin{proof}[1 implies 2]

Define 
\[F(\omega, t, Y_t, Z_t) := Y_t- \phi_{Z_t}^{-1}(Y_t).\]
It follows that 
\[Y_t\mapsto Y_t - F(\omega, t, Y_t, Z_t) = \phi_{Z_t}^{-1}(Y_t)\]
is a bijection for all $Z_t$, and that, as $\phi_{Z_t}$ is is invariant under equivalence $\sim_{M_{t+1}}$, so is $F$.

Therefore, there is a solution to 
\[Y_{t+1} = Y_t - F(\omega, t, Y_t, Z_t) + Z_t M_{t+1},\]
 and this satisfies
\[Y_{t+1}-Z_tM_{t+1} = \phi_{Z_t}^{-1}(Y_t),\]
hence
\[\phi(Y_{t+1}) = \phi_{Z_t}(Y_{t+1}-Z_tM_{t+1}) = Y_t\]
as desired.
\end{proof}
\begin{proof}[2 implies 1]
We know $F$ satisfies the conditions of Theorem \ref{thm:BSDEExist}. We define 
\[\phi(Y_{t+1}) := Y_{t},\]
where $(Y_t, Z_t)$ is the solution at time $t$ of the BSDE (\ref{eq:BSDEgen}) with terminal condition $Y_{t+1}$ at time $t$.

We know this solution is unique. As $F$ is invariant under equivalence $\sim_{M_{t+1}}$, so is $Y_t$ and hence $\phi$. For a fixed $Z_t$, as 
\[Y_{t} \mapsto Y_t-F(\omega, t, Y_t, Z_t) = E[Y_{t+1}|\mathcal{F}_t] = Y_{t+1} - Z_tM_{t+1}\]
is a bijection, clearly $\phi_{Z_t}$ satisfies
\[\phi_{Z_t}(Y_t-F(\omega, t, Y_t, Z_t)) =\phi_{Z_t}(Y_{t+1} - Z_tM_{t+1}) = \phi(Y_{t+1}) = Y_t,\]
that is, $\phi_{Z_t}$ is the inverse of $Y_{t} \mapsto Y_t-F(\omega, t, Y_t, Z_t)$, and is hence also a bijection.
\end{proof}

\begin{remark}
A key point in this Theorem is that the pairs $(Y_t, Y_{t+1})$ can be arbitrary -- any relationship $Y_t=\phi(Y_{t+1})$ is possible, provided the stated properties hold.
\end{remark}

In fact, it is possible to determine the value of $F$ without using $\mathfrak{O}_t$, given slightly different information.

\begin{definition}
Let $F$ be a driver for a discrete time BSDE (\ref{eq:BSDEgen}) satisfying the conditions of Theorem \ref{thm:BSDEExist}. We define $\mathfrak{E}_t$, the time $t$ ``endpoints'' under $F$,  to be the set of ordered pairs 
\[\mathfrak{E}_t:=\{(Y_t, Q)|Q\in L^\infty(\mathcal{F}_T)\},\] where $Y_t$ is the time $t$ solution to the BSDE with driver $F$ and terminal condition $Q$.
\end{definition}

Again, these endpoints are not arbitrary, but are known to come from a BSDE, and therefore satisfy certain consistency properties.

\begin{definition}
Let $F$ be the driver for a discrete time BSDE (\ref{eq:BSDEgen}) satisfying the conditions of Theorem \ref{thm:BSDEExist}. We define $\mathfrak{F}$, the zero hedging function associated with $F$, to be 
\[\mathfrak{F}(\omega, t, Y_t) := F(\omega, t, Y_t, 0).\]
\end{definition}

\begin{remark}
It will become clear that the choice of $Z=0$ in the definition of $\mathfrak{F}$ is arbitrary, and similar results could be obtained given the definition $\mathfrak{F}(\omega, t, Y_t) := F(\omega, t, Y_t, z_t)$
for a given process $z$.
\end{remark}

\begin{definition}
A pair $(\mathfrak{F}, \mathfrak{E}_t)$ is called consistent if, for any initial value $Y_t\in L^\infty(\mathcal{F}_t)$, the the pair $(Y_t, Y_T)\in\mathfrak{E}_t$, where $Y_T$ is the solution of the recursion
\[Y_{u+1}=Y_u+\mathfrak{F}(\omega, u, Y_u).\]
\end{definition}

\begin{lemma}
For $T=1$, the set of endpoints $\mathfrak{E}_0$ will be associated with a unique driver $F$, and hence with a unique consistent zero hedging function $\mathfrak{F}$.
\end{lemma}
\begin{proof}
This is because, in this case, $\mathfrak{E}_0=\mathfrak{O}_0$, and Theorem \ref{thm:OdetF} gives the result.
\end{proof}

\begin{lemma}
For $T\geq 2$, a set of endpoints $\mathfrak{E}_t$ for a given $t\in\{0,...,T-1\}$ will have infinitely many consistent zero hedging functions $\mathfrak{F}$ associated with it.
\end{lemma}
\begin{proof}

Let $\mathfrak{F}$ be a consistent zero hedging function, which exists as $\mathfrak{E}$ comes from some BSDE. For any $k\in\mathbb{R}^K$, define
\[\begin{split}
\tilde{\mathfrak{F}}(\omega, t, Y_t) &= \mathfrak{F}(\omega, t, Y_t)+k\\
\tilde{\mathfrak{F}}(\omega, t+1, Y_{t+1}) &= \mathfrak{F}(\omega, t+1, Y_{t+1}+k)-k,
\end{split}\]
and $\tilde{\mathfrak{F}}(\omega, u, Y_u) = \mathfrak{F}(\omega, u, Y_u)$ for $u>t+1$. Then it is easy to see that $\tilde{\mathfrak{F}}$ is also consistent for $\mathfrak{E}$.
\end{proof}

\begin{theorem} \label{thm:ObsDriver}
Given a consistent pair $(\mathfrak{F}, \mathfrak{E}_t)$, there is a unique function $F(\omega, t, \cdot, \cdot)$ associated with $(\mathfrak{F}, \mathfrak{E}_t)$, that is, where $F(\omega, t, Y_t, Z_t)$ is the value at time $t$, for a given $Y_t$ and $Z_t$, of the driver $F$ generating $\mathfrak{F}$ and $\mathfrak{E}_t$.
\end{theorem}

\begin{proof}
For any $Y_{t+1}\in L^\infty(\mathcal{F}_{t+1})$, we can define a value $Y_T$ through the recursion 
\[Y_{u+1}= Y_u-\mathfrak{F}(\omega, u, Y_u).\]
This value $Y_T$ will appear once and only once in $\mathfrak{E}_t$, in a pair
\[(Y_t, Y_T)\in\mathfrak{E}_t.\]

It is then clear that $Y_t$ and $Y_{t+1}$ are the values (at times $t$ and $t+1$ respectively) of the solution to the BSDE with terminal value $Y_T$, and that the solution $Z$ process will have
\[Z_tM_{t+1} = Y_{t+1}-E[Y_{t+1}|\mathcal{F}_t]\]
and $Z_u=0$ for $u>t$.
We know from Lemma \ref{lem:bijective} that there is a $Y_{t+1}$ associated with every pair $(Y_t, Z_t)$, and therefore, that this will give us all possible pairs $(Y_t, Y_{t+1})$.

We have therefore constructed the set $\mathfrak{O}_t$. The result follows by Theorem \ref{thm:OdetF}.
\end{proof}
\section{Applications to Risk Measures}

We now focus our attention on the theory of risk measures, as in \cite{Follmer2002}, and nonlinear expectations, as in \cite{Peng2005}. The connection between these is present in \cite{Barrieu2005}, and more generally in \cite{Rosazza2006}. We shall not discuss in detail the more general theory of nonlinear evaluations, as in \cite{Peng2005} or \cite{Cohen2008b}. Many of the results in this section parallel those in \cite{Coquet2002}, which discusses continuous time processes related to Brownian motion. 

As in \cite{Cohen2008b}, we follow \cite{Peng2005} by giving the following definition.

\begin{definition} \label{def:NonlinExp}
Let $\{\mathcal{Q}_t\}$ be a family of subsets $\{\mathcal{Q}_t\subset L^2(\mathcal{F}_T)\}$. A system of operators 
\[\mathcal{E}(\cdot|\mathcal{F}_t):L^2(\mathcal{F}_T)\rightarrow L^2(\mathcal{F}_t), \  0\leq t \leq T\]
an $\mathcal{F}_t$-consistent \textbf{nonlinear expectation} for $\{\mathcal{Q}_t\}$ if it satisfies the following properties:
\begin{enumerate}
	\item For $Q, Q'\in\mathcal{Q}_t$, if $Q\geq Q'$ $\mathbb{P}$-a.s. componentwise, then
	\[\mathcal{E}(Q|\mathcal{F}_t) \geq \mathcal{E}(Q'|\mathcal{F}_t)\] $\mathbb{P}$-a.s. componentwise, with, for each $i$,
	\[e_i^*\mathcal{E}(Q|\mathcal{F}_t) = e_i^*\mathcal{E}(Q'|\mathcal{F}_t)\]
	 only if $e_i^*Q = e_i^*Q'$ $\mathbb{P}$-a.s.
	\item $\mathcal{E}(Q|\mathcal{F}_t) = Q$ $\mathbb{P}$-a.s. for any $\mathcal{F}_t$-measurable $Q$.
	\item $\mathcal{E}(\mathcal{E}(Q|\mathcal{F}_t)|\mathcal{F}_s) = \mathcal{E}(Q|\mathcal{F}_s)$ $\mathbb{P}$-a.s. for any $s\leq t$
	\item For any $A\in \mathcal{F}_t$, $I_A\mathcal{E}(Q|\mathcal{F}_t) = \mathcal{E}(I_AQ|\mathcal{F}_t)$ $\mathbb{P}$-a.s.
\end{enumerate}
\end{definition}

We know from \cite{Rosazza2006} that the theory of dynamic risk measures is closely related to the theory of nonlinear expectations. In particular, we can define a dynamic risk measure as the function
\[\rho_t(Q) := -\mathcal{E}(Q|\mathcal{F}_t),\]
where $\mathcal{E}(\cdot|\mathcal{F}_t)$ is an $\mathcal{F}_t$-consistent nonlinear expectation.

Similar ideas are used in \cite{Kloppel2007}, where $\mathcal{E}$ is assumed to be concave, and is referred to as a `monetary concave utility functional'. In the same discrete time, finite state context as considered here, \cite{Jobert2008} considers this structure under the name of a `concave valuation operator'. As shown in \cite{Coquet2002}, in continuous time with Brownian motions, the theory of Backward Stochastic Differential Equations is an appropriate context to study these operators in general.

It is straightforward to prove various properties of risk measures from their definitions in terms of nonlinear expectations. The interested reader is referred to Section 9 of \cite{Cohen2008b} -- the proofs in this discrete time context follow without changes.

One important contribution to the theory of risk measures developed in this paper is that the quantities may be vector valued. Our proofs all work in a multidimensional context, which may be significant in applications related to multiobjective optimisation.

\begin{definition}
 A family of maps $\mathcal{E}(\cdot|\mathcal{F}_t):L^\infty(\mathcal{F}_T)\rightarrow L^\infty(\mathcal{F}_t)$ will be called \textbf{dynamically monotone} for $\{\mathcal{Q}_t\}$ if, for all $s\leq t$,
 
\begin{enumerate}[(i)]
\item $\mathcal{E}$ is an $\mathcal{F}_t$-consistent nonlinear expecation for $\{\mathcal{Q}_t\}$,
\item  $\mathcal{Q}_s\subseteq \mathcal{Q}_t$ ($\{\mathcal{Q}_t\}$ is \textbf{nondecreasing} in $t$)
\item $\mathcal{E}(Q|\mathcal{F}_t)\in\mathcal{Q}_s$ for all $Q\in\mathcal{Q}_s$.
\end{enumerate}
\end{definition}

\begin{remark} In many applications, it may be that $\mathcal{Q}_t= L^\infty (\mathcal{F}_T)$ for all $t$. In this case, a $\mathcal{F}_t$-consistent nonlinear expectation for $\{\mathcal{Q}_t\}$ will automatically be dynamically monotone for $\{\mathcal{Q}_t\}$.
\end{remark}

\begin{definition}
An $\mathcal{F}_t$-consistent nonlinear expectation $\mathcal{E}(.|\mathcal{F}_t)$ is said to be \textbf{(dynamically) translation invariant} if for any $Q\in L^2(\mathcal{F}_T)$, any $q \in L^2(\mathcal{F}_t)$,
\[\mathcal{E}(Q+q|\mathcal{F}_t) = \mathcal{E}(Q|\mathcal{F}_t)+q.\]
\end{definition}

As in \cite{Cohen2008b}, we make the following definition, which ensures a comparison theorem will hold on $[t,T]$ under certain circumstances.
\begin{definition}
Consider some nondecreasing family of sets $\{\mathcal{Q}_t\subset L^2(\mathcal{F}_T)\}$ and some driver $F$, satisfying the assumptions of Theorem \ref{thm:BSDEExist}. Suppose that, for each $t$, for any $Q^1, Q^2\in\mathcal{Q}_t$, the corresponding BSDE solutions $(Y^1,Z^1), (Y^2, Z^2)$ satisfy
\begin{enumerate}[(i)]
\setcounter{enumi}{2}
\item $\mathbb{P}$-a.s., for all $i$, the $i$th component of $F$, given by $e_i^*F$, satisfies
	\[e_i^*F(\omega, t, Y_t^2, Z_t^1) - e_i^*F(\omega, t, Y_t^2, Z_t^2)\geq\min_{j\in \mathbb{J}_t}\{e_i^*(Z^1_t-Z^2_t)(e_j-E[X_{t+1}|\mathcal{F}_{t}])\},\]
	with equality only if $e_i^*Z_t^1\sim_{M_{t+1}}e_i^*Z_t^2$.
	\item $\mathbb{P}$-a.s., if
	\[Y^1_t -F(\omega, t, Y_t^1, Z_t^1) \geq Y^2_t-F(\omega, t, Y_t^2, Z_t^1)\]
	then $Y^1_t\geq Y^2_t$,
	\end{enumerate}
	(cf. assumptions (iii) and (iv) of Theorem \ref{thm:CompThm} and the assumption of Corollary \ref{cor:StrictComp}).
Then we shall call $F$ a \textbf{balanced driver} on $\{\mathcal{Q}_t\}$.
\end{definition}

\begin{remark} \label{rem:timpliesspostt}
Note that as $\{\mathcal{Q}_t\}$ is nondecreasing, this implies that, for any $Q^1, Q^2\in\mathcal{Q}_t$, for all $i$, the corresponding BSDE solutions $(Y^1,Z^1), (Y^2, Z^2)$ satisfy
	\[e_i^*F(\omega, s, Y_s^2, Z_s^1) - e_i^*F(\omega, s, Y_s^2, Z_s^2)\geq\min_{j\in \mathbb{J}_s}\{e_i^*(Z^1_s-Z^2_s)(e_j-E[X_{s+1}|\mathcal{F}_{s}])\},\]
	with equality only if $e_i^*Z_s^1\sim_{M_{s+1}}e_i^*Z_s^2$, for all $s>t$. Similarly for Assumption (iv).
	\end{remark}

\begin{remark}
As in Corollary \ref{cor:CompThmAssnSwitch}, these assumptions can be generalised slightly by changing which of $Y^1$ and $Y^2$, and which of $Z^1$ and $Z^2$, appears in each place.
\end{remark}

In this context of a discrete time space generated by a finite state system, we establish the following theorem, which directly relates BSDEs (and their drivers) to nonlinear expectations. This theorem applies in both the scalar and vector cases.

\begin{theorem} \label{thm:EqivBSDENonlin}
For some family of operators $\mathcal{E}(.|\mathcal{F}_t)$, let $\{\mathcal{Q}_t\subset L^2(\mathcal{F}_T)\}$,  be such that if $Q\in\mathcal{Q}_t$ then $Q+q\in \mathcal{Q}_t$ for all $q\in L^\infty(\mathcal{F}_t)$, and $\mathcal{E}$ is dynamically monotone for $\{\mathcal{Q}_t\}$. Then the following two statements are equivalent.
 \begin{enumerate}
 \item $\mathcal{E}(.|\mathcal{F}_t)$ is an $\mathcal{F}_t$-consistent, dynamically translation invariant, nonlinear expectation for $\{\mathcal{Q}_t\}$.
 \item There exists a driver $F$, which is balanced on $\{\mathcal{Q}_t\}$, is independent of $Y$, and satisfies the normalisation condition $F(\omega, t, Y_t, 0) = 0$, such that, for all $Q$, $Y_t = \mathcal{E}(Q|\mathcal{F}_t)$ is the solution to a BSDE with terminal condition $Q$ and driver $F$.
 \end{enumerate}
 Furthermore, these two statements are related by the equation
 \[F(\omega, t, Y_{t}, Z_{t}) = \mathcal{E}(Z_{t}M_{t+1}|\mathcal{F}_t).\]
 \end{theorem}

\begin{proof}[2 implies 1.]
Let $\mathcal{E}(Q|\mathcal{F}_t):=Y_t$, the time $t$ solution of the BSDE with terminal value $Q$. We shall show that each of the properties of a nonlinear expectation is satisfied.
\begin{enumerate} 
	\item The statement $\mathcal{E}(Q^1|\mathcal{F}_t) \geq \mathcal{E}(Q^2|\mathcal{F}_t)$ $\mathbb{P}$-a.s. whenever $Q^1, Q^2\in \mathcal{Q}$, $Q^1\geq Q^2$ $\mathbb{P}$-a.s. is the main result of Theorem \ref{thm:CompThm}, which holds on $[t,T]$ as $F$ is balanced on $\mathcal{Q}_t$ (see Remark \ref{rem:timpliesspostt}). The strict comparison of Corollary \ref{cor:StrictComp} then establishes the second statement.
	\item By normalisation, the solution to the BSDE with $\mathcal{F}_t$-measurable terminal condition $Q$ will be $(Y_s, Z_s) = (Q, 0)$ for $s\geq t$. By the uniqueness result of Theorem \ref{thm:BSDEExist} this is then the value of $\mathcal{E}(Q|\mathcal{F}_t)$.
	\item As $Y$ is the solution to the relevant BSDE, we can deduce
  \[Y_t= Y_s - \sum_{s\leq u< t} F(\omega, u, Y_u, Z_u) + \sum_{s\leq u< t} Z_u M_{u+1}.\]
  Hence $Y_s$ is also the time $s$ value of a solution to the BSDE with terminal time $t$ and value $Y_t$. Hence 
  $\mathcal{E}(\mathcal{E}(Q|\mathcal{F}_t)|\mathcal{F}_s) = \mathcal{E}(Q|\mathcal{F}_s)$ $\mathbb{P}$-a.s. as desired.
	\item We know that 
	\[I_AQ = I_A Y_t - \sum_{t\leq u< T} I_AF(\omega, u, Y_u, Z_u) + \sum_{t\leq u< T} I_AZ_u M_{u+1}\]
	and by normalisation, 
	\[I_AF(\omega, u, Y_u, Z_u)= F(\omega, u, I_AY_u, I_AZ_u).\]
	
	Hence $(I_A Y, I_A Z)$ is the solution to a BSDE with driver $F$ and terminal condition $I_A Q$, and hence
	\[ I_A\mathcal{E}(Q|\mathcal{F}_t) =I_A Y_t = \mathcal{E}(I_AQ|\mathcal{F}_t)\]
	as desired.
\end{enumerate}
	
We can also show that this nonlinear expectation is dynamically translation invariant. For any $Q$ let $(Y, Z)$ be the solution to
\[Q = Y_t - \sum_{t\leq u< T} F(\omega, u, Y_u, Z_u) + \sum_{t\leq u< T} Z_u M_{u+1}.\]
Then, for the terminal condition $Q+q$,
\[Q+q = Y_t+q - \sum_{t\leq u< T} F(\omega, u, Y_u+q, Z_u) + \sum_{t\leq u< T} Z_u M_{u+1},\]
 as $F$ is independent of $Y$. Hence we have that $(Y+q, Z)$ is the solution on $[t,T]$ for the BSDE with terminal condition $Q+q$. The result follows.
 
Finally, we can see that 
\[Y_{t+1}= Y_t - F(\omega, t, Y_t, Z_t) + Z_tM_{t+1}.\]
Taking an $\mathcal{F}_t$-conditional expectation and rearranging gives
\[\begin{split}
F(\omega, t, Y_t, Z_t)&= Y_t -E[Y_{t+1}|\mathcal{F}_t]\\
&=\mathcal{E}(Y_{t+1}|\mathcal{F}_t) -E[Y_{t+1}|\mathcal{F}_t]\\
&=\mathcal{E}(Y_{t+1}-E[Y_{t+1}|\mathcal{F}_t]|\mathcal{F}_t) \\
&=\mathcal{E}(Z_{t}M_{t+1}|\mathcal{F}_t) 
\end{split}\]
as desired.
\end{proof}

\begin{proof}[1 implies 2.]
We know that, for any $0\leq t<T$, we can write \[\mathcal{E}(Q|\mathcal{F}_t) = \mathcal{E}(\mathcal{E}(Q|\mathcal{F}_{t+1})|\mathcal{F}_t).\]

We propose that $\mathcal{E}(Q|\mathcal{F}_t)$ will satisfy a BSDE with driver \[F(\omega, t, Y_t, Z_{t}) := \mathcal{E}(Z_tM_{t+1}|\mathcal{F}_t).\] 

First, we note that for any $t$ and any $Q$,
\[\mathcal{E}(Q|\mathcal{F}_{t+1}) - E[\mathcal{E}(Q|\mathcal{F}_{t+1})|\mathcal{F}_t]\]
is an $\mathcal{F}_{t+1}$-measurable random variable with $\mathcal{F}_t$-conditional mean zero. Therefore we can apply Theorem \ref{thm:MartRep} to show that, for some $\mathcal{F}_t$-measurable matrix $Z_t$,
\begin{equation}\label{eq:p21e1}
\mathcal{E}(Q|\mathcal{F}_{t+1}) - E[\mathcal{E}(Q|\mathcal{F}_{t+1})|\mathcal{F}_t]=Z_tM_{t+1}.
\end{equation}

As $\mathcal{E}$ is dynamically translation invariant,
\begin{equation}\label{eq:p21e2}
\begin{split}
F(\omega, t, Y_t, Z_{t}) &= \mathcal{E}(Z_tM_{t+1}|\mathcal{F}_t)\\
&=\mathcal{E}(\mathcal{E}(Q|\mathcal{F}_{t+1})|\mathcal{F}_t) - E[\mathcal{E}(Q|\mathcal{F}_{t+1})|\mathcal{F}_t]\\
&=\mathcal{E}(Q|\mathcal{F}_t) - E[\mathcal{E}(Q|\mathcal{F}_{t+1})|\mathcal{F}_t].
\end{split}\end{equation}

Therefore, if we define $Y_t:=\mathcal{E}(Q|\mathcal{F}_t)$, we can combine (\ref{eq:p21e1}) and (\ref{eq:p21e2}) to give
\[Y_{t+1} = Y_t - F(\omega, t, Y_t, Z_{t}) + Z_tM_{t+1}.\]

The one-step dynamics being established, the $\mathcal{E}(Q|\mathcal{F}_t)$ satisfies the BSDE with driver $F$ by induction. 

We need only to show that the driver $F$ is balanced on $\mathcal{Q}_t$ for all $t$. As $F$ is independent of $Y$, the only relevant requirement for $F$ to be balanced is that, for each component $i$, for $Z_t^1, Z^2_t$ defined as above for $Q^1, Q^2\in\mathcal{Q}_t$, 
\[e_i^*F(\omega, t, Y_t^2, Z_t^1) - e_i^*F(\omega, t, Y_t^2, Z_t^2)\geq\min_{j\in \mathbb{J}_t}\{e_i^*(Z^1_t-Z^2_t)(e_j-E[X_{t+1}|\mathcal{F}_{t}])\},\]
	with equality only if $e_i^*Z_t^1\sim_{M_{t+1}}e_i^*Z_t^2$.
	
Let $Y^k_t = \mathcal{E}(Q^k|\mathcal{F}_t)$ for $k=1,2$. Define an $\mathcal{F}_t$-measurable random variable $q$ by

\[\begin{split}
e_i^*q &= \min_{j\in \mathbb{J}_t}\{e_i^*(Y^1_{t+1}-Y^2_{t+1})|\mathcal{F}_t, X_{t+1} = e_j\}\\
&=E[Y^1_{t+1}-Y^2_{t+1}|\mathcal{F}_t]+\min_{j\in \mathbb{J}_t}\{e_i^*(Z^1_t-Z^2_t)(e_j-E[X_{t+1}|\mathcal{F}_{t}])\}.\end{split}
\]

Then $Y^1_{t+1}-q\geq Y^2_{t+1}$ componentwise, and hence, as these values lie in $\mathcal{Q}_t$, (by dynamic monotonicity and the fact $q$ is $\mathcal{F}_t$-measurable), we know 
\[\mathcal{E}(Y^1_{t+1}-q|\mathcal{F}_t)\geq \mathcal{E}(Y^2_{t+1}|\mathcal{F}_t)\]

By dynamic translation invariance of $\mathcal{E}$, we then have, 
\[e_i^*\mathcal{E}(Y^1_{t+1}|\mathcal{F}_t)- e_i^*\mathcal{E}(Y^2_{t+1}|\mathcal{F}_t)\geq e_i^*q\]
and hence, subtracting $e_i^*E[Y^1_{t+1}-Y^2_{t+1}|\mathcal{F}_t]$ fom both sides,
\[e_i^*\mathcal{E}(Z^1_tM_{t+1}|\mathcal{F}_t)- e_i^*\mathcal{E}(Z^2_tM_{t+1}|\mathcal{F}_t)\geq \min_{j\in \mathbb{J}_t}\{e_i^*(Z^1_t-Z^2_t)(e_j-E[X_{t+1}|\mathcal{F}_{t}])\}.\]
Given $F(\omega, t, Y_t, Z_t) = \mathcal{E}(Z_tM_{t+1}|\mathcal{F}_t)$, this shows
\[e_i^*F(\omega, t, Y^1_t, Z^1_t)- e_i^*F(\omega, t, Y^2_t, Z^2_t)\geq \min_{j\in \mathbb{J}_t}\{e_i^*(Z^1_t-Z^2_t)(e_j-E[X_{t+1}|\mathcal{F}_{t}])\}\]
as desired.

To show the strict inequality, note that, for each $i$, if 
\[e_i^*\mathcal{E}(Y^1_{t+1}-q|\mathcal{F}_t) =e_i^*\mathcal{E}(Y^2_{t+1}|\mathcal{F}_t)\]
then $e_i^*Y^1_{t+1}-e_i^*q=e_i^*Y^2_{t+1}$ $\mathbb{P}$-a.s. It follows that 
\[e_i^*(Y^1_{t+1}-Y^2_{t+1}) = e_i^*q = \min_{j\in \mathbb{J}_t}\{e_i^*(Y^1_{t+1}-Y^2_{t+1})|\mathcal{F}_t, X_{t+1} = e_j\}\]
and so $e_i^*(Y^1_{t+1}-Y^2_{t+1})=0$ $\mathbb{P}$-a.s. It follows that $e_i^*Z^1_t\sim_{M_{t+1}} e_i^*Z^2_t$ as required.
\end{proof}

Theorem \ref{thm:EqivBSDENonlin} is significant in that, for both the scalar and vector cases, it characterises all dynamically translation invariant $\mathcal{F}_t$-consistent nonlinear expectations as the solutions to BSDEs.  These are the nonlinear expectations which are associated with dynamically consistent risk measures. Unlike the continuous time context with Brownian motions considered in \cite{Coquet2002}, no further assumptions on $F$ or $\mathcal{E}$ are needed in this discrete context.

\section{From Static to Dynamic behaviour} \label{sec:StaticDynamic}
A final question that remains to be answered is whether, given a particular map $Q\mapsto Y_0$, it is possible to find an $\mathcal{F}_t$-consistent nonlinear expectation which agrees with it. The following theorem addresses this problem in the scalar case, under the assumption of monotonicity.

\begin{theorem} \label{thm:dynamicism}
Consider a measurable, scalar valued map $\mathcal{E}:L^2(\mathcal{F}_T) \rightarrow L^2(\mathcal{F}_0)$. Suppose this map satisfies:
\begin{enumerate}[(i)]
\item ($\mathcal{F}_t$-consistency) For any $Q\in L^2(\mathcal{F}_T)$ and any $t\leq T$, there exists an $\mathcal{F}_t$-measurable random variable $Y_t$ such that 
\begin{equation} \label{eq:FconsistY}
\mathcal{E}(I_A Q) = \mathcal{E}(I_A Y_t)
\end{equation}
for any $A\in \mathcal{F}_t$.
\item ($\mathcal{F}_0$-triviality) $\mathcal{E}(Q) = Q$ for all $\mathcal{F}_0$-measurable $Q$.
\item (Monotonicity) For any $Q, Q'\in L^2(\mathcal{F}_t)$, if $Q\geq Q'$ $\mathbb{P}$-a.s., then $\mathcal{E}(Q) \geq \mathcal{E}(Q')$, with equality if and only if $Q=Q'$ $\mathbb{P}$-a.s.
\end{enumerate}
Then there exists a unique  
dynamic nonlinear expectation $\mathcal{E}(.|\mathcal{F}_t)$ (for $\mathcal{Q}=L^2(\mathcal{F}_T)$) such that 
\[\mathcal{E}(Q)=\mathcal{E}(Q|\mathcal{F}_0)\]
for all $Q\in L^2(\mathcal{F}_T)$. This dynamic nonlinear expectation is given by 
\[\mathcal{E}(Q|\mathcal{F}_t) = Y_t\]
with $Y_t$ as in $(i)$.
\end{theorem}

\begin{proof}
We shall first show that the random variable $Y_t$ satisfying $\mathcal{F}_t$-consistency is unique. We know from the assumption of $\mathcal{F}_t$-consistency that some $Y_t$ exists satisfying (\ref{eq:FconsistY}).

As $\mathcal{F}$ is generated by $X$, we can consider the atomic event $A$ given by
\[A=\{X_0=e_{i_0}, X_1=e_{i_1} ..., X_t=e_{i_t}\}\]
for the path $(e_{i_0},..., e_{i_t})$, where $e_{i_s}$ is a basis vector in $\mathbb{R}^N$, for all $s$. As $A$ is an atomic event in $\mathcal{F}_t$, the $\mathcal{F}_t$-measurable random variable $Y_t$ must be constant on $A$. Let $y_A\in\mathbb{R}$ be the value $Y_t$ takes on $A$, i.e. $I_Ay_A=I_AY_t$ $\mathbb{P}$-a.s.

Hence if $y=y_A$,
\begin{equation}\label{eq:atomicQ}
\mathcal{E}(I_A Q) = \mathcal{E}(I_A Y_t) = \mathcal{E}(I_A y).
\end{equation}
By monotonicity, the map 
\[y \mapsto \mathcal{E}(I_A y)\]
is strictly increasing as a function of $y$, and, therefore, there is a unique value $y$ which solves (\ref{eq:atomicQ}). Clearly, this is given by $y=y_A$. Hence, if $Y^1_t$ and $Y^2_t$ both satisfy (\ref{eq:FconsistY}), then $I_AY^1_t = I_AY^2_t = I_Ay_A$ $\mathbb{P}$-a.s. As $\mathcal{F}_t$ is generated by a finite number of events of the form of $A$, it follows that $Y^1_t=Y^2_t = \sum_A I_A y_A$ $\mathbb{P}$-a.s. Hence the random variable $Y_t$ satisfying $\mathcal{F}_t$-consistency is unique (up to equality $\mathbb{P}$-a.s.).

If $A=\Omega$, then the assumption of $\mathcal{F}_t$-consistency at $t=0$, along with $\mathcal{F}_0$-triviality, implies that 
\[\mathcal{E}(Q) = \mathcal{E}(Y_0) = Y_0\]
and therefore if $\mathcal{E}(Q|\mathcal{F}_t):=Y_t$, we have
\[\mathcal{E}(Q)= \mathcal{E}(Q|\mathcal{F}_0).\]

We now wish to show that $Y_t$ satisfies the properties of a nonlinear expectation, as given in Definition \ref{def:NonlinExp}.
\begin{enumerate}
\item Suppose $Q\geq Q'$ $\mathbb{P}$-a.s. Let $Y$ and $Y'$ be the corresponding processes from the $\mathcal{F}$ consistency assumption. Then we know from monotonicity that, for any $t\leq T$, $\mathcal{E}(I_AQ)\geq \mathcal{E}(I_AQ')$ for all $A\in\mathcal{F}_t$. By the same argument as used to show uniqueness, this implies that $I_AY_t\geq I_AY_t'$ for all atomic $A$, and hence $Y_t\geq Y_t'$ $\mathbb{P}$-a.s. Hence 
\[Y_t=\mathcal{E}(Q|\mathcal{F}_t)\geq \mathcal{E}(Q'|\mathcal{F}_t)=Y_t'\]
$\mathbb{P}$-a.s. as desired.
\item If $Q$ is $\mathcal{F}_t$-measurable, then $Y_s=Q$ for all $t\leq s\leq T$ satisfies the assumptions of $\mathcal{F}_t$-consistency. We have established that this solution is unique, and hence
\[\mathcal{E}(Q|\mathcal{F}_t) = Y_t = Q\]
$\mathbb{P}$-a.s. as desired.
\item Let $s\leq t$. For any $Q\in L^2(\mathcal{F}_T)$, let $Y$ be the corresponding process for $\mathcal{E}(Q|\mathcal{F}_{(.)})$. Let $\tilde Y$ be the process for $\mathcal{E}(Y_t|\mathcal{F}_{(.)})$. Then we know that $\tilde Y$ satisfies
\[\mathcal{E}(I_A \tilde Y_s) = \mathcal{E}(I_A Y_t) = \mathcal{E}(I_A Q)\]
for all $A\in\mathcal{F}_s$. Therefore $\tilde Y_s$ is also a solution for $\mathcal{E}(Q|\mathcal{F}_s)$. By uniqueness, this implies that $\tilde Y_s = Y_s$ $\mathbb{P}$-a.s. Therefore, by definition,
\[\mathcal{E}(\mathcal{E}(Q|\mathcal{F}_t)|\mathcal{F}_s)=\mathcal{E}(Y_t|\mathcal{F}_s) = \tilde Y_s= Y_s = \mathcal{E}(Q|\mathcal{F}_s)\]
$\mathbb{P}$-a.s. as desired.
\item Fix $t\leq T$. We need to show that, for any $A\in\mathcal{F}_t$, if $Y$ is the process associated with $Q$ and $\tilde Y$ is the process associated with $I_A Q$, then 
\[I_A Y_t = \tilde Y_t.\] 

For any $B\in\mathcal{F}_t$, $\mathcal{F}_t$-consistency shows that 
\[\mathcal{E}(I_B I_{A} Y_t) =\mathcal{E}(I_{A\cap B} Y_t)=\mathcal{E}(I_{A\cap B} Q)= \mathcal{E}(I_B I_A Q) = \mathcal{E}(I_B \tilde Y_t).\]
Therefore, $I_{A} Y_t$ is also a solution to the $\mathcal{F}_t$-consistency assumption for $I_A Q$. By uniqueness, this shows that $I_A Y_t = \tilde Y_t$ $\mathbb{P}$-a.s. as desired.
\end{enumerate}

Therefore, we have shown that $\mathcal{E}(\cdot|\mathcal{F}_t)=Y_t$ is a dynamic nonlinear expectation for $\mathcal{Q}=L^2(\mathcal{F}_T)$.
\end{proof}
\begin{corollary} \label{cor:DynamCondNec}
The assumptions of Theorem \ref{thm:dynamicism} are necessary.
\end{corollary}
\begin{proof}
The assumptions of Theorem \ref{thm:dynamicism} are all properties of a dynamic nonlinear expectation (for $\mathcal{Q}=L^2(\mathcal{F}_T)$). $\mathcal{F}_t$-consistency is a consequence of a combination of Properties 3 and 4 of nonlinear expectations, $\mathcal{F}_0$-triviality is simply an application of Property 2 at $t=0$ and monotonicity is simply an application of Property 1 at $t=0$.
\end{proof}

\begin{remark}
Unfortunately, it is not true, in general, that $\mathcal{F}_t$-consistency is satisfied for an arbitrary nonlinear expectation. A counterexample, of an $\mathcal{F}_t$-inconsistent nonlinear expecation, is given below (Example \ref{ex:Counter}). In the case where $\mathcal{E}=E_{\tilde{\mathbb{P}}}$, the classical expectation under a measure $\tilde{\mathbb{P}}$, this statement can be shown to hold using the Radon-Nikodym theorem, however this does not extend in a straightforward manner to the nonlinear problems considered here.
\end{remark}

\begin{example}\label{ex:Counter}
Consider the simple case where $X$ can take two values at each time point with equal probability, and $T=2$. Consider the nonlinear expecation given by 
\[\mathcal{E}(X) = 0.1\times E(X) + 0.9 \times \inf_{x\in\mathbb{R}}\{x|\mathbb{P}(X_T\leq x)>0\}.\]
We shall show that this nonlinear expectation is not $\mathcal{F}_t$ consistent, and, therefore, there is no $\mathcal{F}_t$-consistent nonlinear expectation (as in Definition \ref{def:NonlinExp}) which agrees with it.

Consider a terminal condition $Q$ with values
\[Q=\left\{\begin{array}{cl} 
0& \text{if }X_1=e_1, X_2=e_1\\
-2& \text{if }X_1=e_1, X_2=e_1\\
4& \text{if }X_1=e_2, X_2=e_1\\
-1& \text{if }X_1=e_2, X_2=e_2.\end{array}\right.\]

We can then see that
\[\begin{split}
\mathcal{E}(I_{X_1=e_1}Q) &= -1.85,\\
\mathcal{E}(I_{X_1=e_2}Q) &= -0.825
\end{split}
\]
and
\[\mathcal{E}(Q)=-1.775.\]

We wish to find an $\mathcal{F}_1$-measurable $Y$ satisfying the requirements of $\mathcal{F}$-consistency. Solving numerically, (values to four decimal places),
\[\mathcal{E}(I_{X_1=e_1}Q)=\mathcal{E}(I_{X_1=e_1}Y) \text{ implies } Y=-2.0556 \text{ given }X_1=e_1.\] 
\[\mathcal{E}(I_{X_1=e_2}Q)=\mathcal{E}(I_{X_1=e_1}Y) \text{ implies } Y=-0.9167 \text{ given }X_1=e_2.\]
These solutions are unique by the monotonicity of $\mathcal{E}$.

However, for these values of $Y$, $\mathcal{E}(Y) = -1.9417 \neq -1.775$. Therefore, there is no $Y$ satisfying the requirements of $\mathcal{F}_t$-consistency for this value of $Q$. Hence $\mathcal{E}$ is not $\mathcal{F}_t$-consistent. 

It follows, from Corollary \ref{cor:DynamCondNec}, that there is no $\mathcal{F}_t$-consistent nonlinear expectation which agrees with $\mathcal{E}$ at time $t=0$.
\end{example}

\bibliographystyle{plain}  
\bibliography{General}
\end{document}